\documentclass{amsart}

\numberwithin{equation}{section}

\usepackage{hyperref}
\usepackage{comment}

\newtheorem{thm}{Theorem}[section]
\newtheorem{lem}[thm]{Lemma}
\newtheorem{cor}[thm]{Corollary}

\newtheorem{rem}[thm]{Remark}

\newcommand\D{{\mathbb D}}

\newcommand\C{{\mathbb C}}

\newcommand\N{{\mathbb N}}

\newcommand\R{{\mathbb R}}

\renewcommand{\epsilon}{\varepsilon}
\renewcommand{\phi}{\varphi}

\DeclareMathOperator{\Leb}{Leb}
\DeclareMathOperator{\Var}{Var}
\newcommand{\norm}[1]{\left\| #1 \right\|}
\newcommand{\moins}{\backslash}
\DeclareMathOperator{\dLeb}{dLeb}
\renewcommand{\hat}{\widehat}
\newcommand{\dd}{\, {\rm d}}

\newcommand\Lp{{\mathcal{BV}}}
\newcommand\TLp{{\mathcal L}}

\begin{document}

\title{A Borel-Cantelli lemma for intermittent interval maps}
\author{S\'{e}bastien Gou\"{e}zel}
\address{
IRMAR\\
Universit\'{e} de Rennes 1\\
Campus de Beaulieu, b\^{a}timent 22\\
35042 Rennes Cedex, France.}
\email{sebastien.gouezel@univ-rennes1.fr}
\date{March 9, 2007}
 \thanks{We thank
Fran\c{c}ois Maucourant for useful comments.}

\begin{abstract}
We consider intermittent maps $T$ of the interval, with an
absolutely continuous invariant probability measure $\mu$. Kim
showed that there exists a sequence of intervals $A_n$ such
that $\sum \mu(A_n)=\infty$, but $\{A_n\}$ does not satisfy the
dynamical Borel-Cantelli lemma, i.e., for almost every $x$, the
set $\{n : T^n(x)\in A_n\}$ is finite. If $\sum
\Leb(A_n)=\infty$, we prove that $\{A_n\}$ satisfies the
Borel-Cantelli lemma. Our results apply in particular to some
maps $T$ whose correlations are not summable.
\end{abstract}

\maketitle

\section{Introduction}

Let $T$ be an ergodic probability preserving transformation of a
space $(X,\mu)$, and let $A_n$ be a sequence of subsets of $X$ with
$\sum \mu(A_n)=+\infty$. It is an interesting question to know
whether, for almost every point $x$, $T^n(x)$ belongs to $A_n$
infinitely often. By the classical Borel-Cantelli lemma, this holds
if the sets $T^{-n}A_n$ are pairwise independent, but this condition
is almost never satisfied for dynamical systems, so one is led to
looking for weaker conditions.

If $T$ is invertible, taking $A_n=T^n(A)$ for some fixed set
$A$ gives a trivial counterexample (and similar counterexamples
also exist for noninvertible maps) . Hence, some regularity
conditions on the sets $A_n$ are necessary. For uniformly
hyperbolic dynamical systems, Chernov and Kleinbock have solved
the problem for lots of families of balls in
\cite{ChernovKleinbock} (see also
\cite{maucourant:BorelCantelli}). The partially hyperbolic case
is dealt with in \cite{dolgopyat:limit}. Concerning
non-uniformly hyperbolic (or expanding) systems, Kim has
considered in \cite{kim:BorelCantelli} a family of interval
maps with a neutral fixed points and obtained partial results.
Our goal in this note is to complete these results (for the
same family of maps) and obtain a full description of the
situation.

Consider some parameter $\alpha>0$ and let $T_\alpha:(0,1]\to(0,1]$
be given by
  \begin{equation}
  T_\alpha(x)=
  \begin{cases}
  x(1+2^\alpha x^\alpha)&\text{ if }x\in (0, 1/2],\\
  2x-1 &\text{ if }x\in (1/2,1].
  \end{cases}
  \end{equation}
It preserves a unique (up to multiplication by a scalar) absolutely
continuous measure $\mu$, and this measure has finite mass if and
only if $\alpha<1$. Henceforth, we will only consider this case, and
assume that $\mu$ is normalized to be a probability measure. We will
also denote by $\Leb$ the Lebesgue measure on $(0,1]$.

In \cite{kim:BorelCantelli}, Kim proves the following result: for
any $\alpha<1$, there exist intervals $A_n$ such that $\sum
\mu(A_n)=\infty$ but, for almost every $x$, $T_\alpha^n(x)\in A_n$
occurs only finitely many times. In other words, the answer to the
Borel-Cantelli problem in this setting is not always positive. On
the other hand, he proves that, if $A_n$ is a sequence of intervals
in $(d,1]$ for some $d>0$, with $\sum \mu(A_n)=\infty$, and
\begin{itemize}
\item either $A_{n+1}\subset A_n$ for all $n$
\item or $\alpha<(3-\sqrt{5})/2$
\end{itemize}
then, for almost every $x$, $T_\alpha^n(x)$ belongs to $A_n$
infinitely many times. In this note, we prove the following
theorem:
\begin{thm}
\label{main_thm} Let $\alpha<1$, and let $A_n$ be a sequence of
intervals with $\sum \Leb(A_n)=\infty$. Then, for almost every $x$,
$T_\alpha^n(x)$ belongs to $A_n$ infinitely many times.
\end{thm}
The measures $\mu$ and $\Leb$ are uniformly equivalent on every
interval $(d,1]$ (more precisely, on every interval $(d,1]$,
the density $h$ of $\mu$ with respect to $\Leb$ is Lipschitz
continuous and bounded from above and below). Hence, this
theorem implies the aforementioned result of Kim.

\smallskip

The proof involves a measurement of how sets $T_\alpha^{-i}A_i$
and $T_\alpha^{-j}A_j$ are ``close to be independent''. For the
following informal description of the proof, assume for the
sake of simplicity that the intervals $A_n$ are all contained
in $(1/2,1]$. The speed of decay of correlations of the map
$T_\alpha$ is exactly $1/n^{\beta-1}$ for $\beta=1/\alpha$,
which means that the best estimate we could hope for is of the
form
  \begin{equation}
  \label{eq:BestDecor}
  |\mu(T_\alpha^{-i}A_i \cap T_\alpha^{-j}A_j) -\mu(A_i)\mu(A_j) |
  \leq \frac{C \mu (A_j)}{(j-i)^{\beta-1}},
  \end{equation}
for $j>i$. This estimate indeed holds, and implies Theorem
\ref{main_thm} when the sequence $1/n^{\beta-1}$ is summable,
that is, when $\alpha<1/2$. However, it is not sufficient when
$1/2\leq \alpha<1$, and we need to know further terms in the
asymptotics of the correlations. Here comes into play our main
technical tool, the renewal sequence of transfer operators,
studied by Sarig in \cite{sarig:decay}. Using the results in
\cite{gouezel:decay}, we will prove the existence of a sequence
$c_n$ converging to $1$ such that
  \begin{equation}
  \label{eq:SubtleBestDecor}
  |\mu(T_\alpha^{-i}A_i \cap T_\alpha^{-j}A_j) -
  c_{j-i}\mu(A_i)\mu(A_j) | \leq \frac{C
  \mu (A_j)}{(j-i)^{\beta}}.
  \end{equation}
This sequence is of the form $c_n=1+c/n^{\beta-1}+o(1/n^{\beta-1})$
for some nonzero constant $c$, which shows that \eqref{eq:BestDecor}
is indeed optimal. For the purposes of the Borel-Cantelli problem,
\eqref{eq:SubtleBestDecor} is sufficient and will imply Theorem
\ref{main_thm} in all cases, since the sequence $1/n^\beta$ is
summable whenever $\alpha<1$.

On the technical level, the results in \cite{sarig:decay,
gouezel:decay} deal with spaces of Lipschitz functions. However, the
essential results are formulated in an abstract Banach spaces
framework. They can therefore also be applied to spaces of functions
with bounded variation, which is what is needed here to deal with
the characteristic functions of intervals.

\begin{rem}
Theorem \ref{main_thm} still holds for transformations with an even
more neutral fixed point, as soon as there is still an absolutely
continuous invariant probability measure. This is for example the
case if the fixed point is of the form $x+x^2 (\log x)^2$, or more
generally for the class of maps introduced by Holland in
\cite{holland}. However, the results of \cite{gouezel:decay} are not
sufficient to prove this, and one needs to use results in the
unpublished thesis \cite{gouezel:these}, for example Remark 2.4.8 or
Remark 2.4.11.
\end{rem}

\section{Abstract tools}

First of all, let us recall a criterion implying the Borel-Cantelli
property (proved e.g.~in \cite[Proposition
6.26.3]{borel_cantelli_dependant}):
\begin{thm}
\label{thm:BC_abstrait}
Let $B_n$ be sets of a probability space
$(X,\mu)$ with $\sum \mu(B_n)=\infty$. Assume that
  \begin{equation}
  \limsup_{n\to \infty}\frac{ \sum_{0\leq i<j<n} \mu(B_i\cap
  B_j)}{\left( \sum_{j=0}^{n-1} \mu(B_j)\right)^2} \leq \frac{1}{2}.
  \end{equation}
Then almost every point of $X$ belongs to infinitely many $B_n$'s.
\end{thm}

We will apply this result to $B_n=T_\alpha^{-n}(A_n)$. Hence,
we need a good quantitative estimate on $\mu(T_\alpha^{-i}A_i
\cap T_\alpha^{-j}A_j)$. This estimate will be provided by
\emph{renewal sequences of transfer operators}, as used by
Sarig in \cite{sarig:decay}. For our purpose, the following
abstract result will be most useful. Let $\overline{\D}$ be the
closed unit disk in $\C$.
\begin{thm}
\label{thm:renewal} Let $\Lp$ be a Banach space, and let
$(R_n)_{n\geq 1}$ be a sequence of continuous linear operators
on $\Lp$. Assume that, for some $\beta>1$,
$\sum_{k>n}\norm{R_k}=O(1/n^\beta)$. Hence, $R(z)=\sum R_n z^n$
and $R'(z)=\sum n R_n z^{n-1}$ are well defined operators on
$\Lp$, for $z\in \overline{\D}$. Assume moreover that $1$ is a
simple isolated eigenvalue of $R(1)$, and that the
corresponding eigenprojector $P$ satisfies $PR'(1)P=\gamma P$
for some $\gamma\not=0$. Assume also that, for any $z\in
\overline{\D} \moins \{1\}$, $I-R(z)$ is invertible on $\Lp$.

Let $T_n= \sum_{l=1}^\infty\sum_{k_1+\dots+k_l=n}R_{k_1}\dots
R_{k_l}$. This operator acts continuously on $\Lp$. Then there
exists a sequence $c_n\in \C$ converging to $1$ such that $T_n
-c_n P/\gamma=O(1/n^\beta)$.
\end{thm}
\begin{proof}
\cite[Theorem 5.4]{gouezel:decay} (for large enough $N$) shows
that $T_n$ converges to $P/\gamma$, and that there exists a
sequence of operators $Q_n$ such that $T_n-PQ_nP=O(1/n^\beta)$.
This theorem even gives a closed form expression for $Q_n$, but
we will not need it.

Since $P$ is a one-dimensional projection, there exists a
complex number $d_n$ such that $PQ_nP=d_nP$. The convergence of
$T_n$ to $P/\gamma$ shows that $d_n$ converges to $1/\gamma$.
We obtain the theorem for $c_n=\gamma d_n$.
\end{proof}

In \cite{sarig:decay, gouezel:decay}, this theorem is applied by
taking $R_n$ to be the ``first return transfer operators'' to
$Y=(1/2,1]$, acting on the space of Lipschitz continuous functions on
$Y$. Here, we will use the same operators $R_n$, but we will use for
$\Lp$ the space of functions of bounded variation on $Y$.

\section{Proof of the main theorem}

In all this section, we fix $\alpha\in (0,1)$ and write $T$ for
$T_\alpha$. Let also $\beta=1/\alpha$.

Let $Y=(1/2,1]$, let $\phi:Y\to\N^*$ denote the first return
time from $Y$ to itself. Let also $\hat{T}$ be the transfer
operator associated to $T$, given for $f\in L^1(\Leb)$ by
  \begin{equation}
  \hat T f(x)=\sum_{T y=x} f(y)/T'(y).
  \end{equation}
Let $R_n f= \hat T^n (1_{Y\cap \{\phi=n\}} f)$, and $T_n f= 1_Y
\hat T^n(1_Y f)$. These operators act on $L^1(Y)$. Moreover,
$R_n$ corresponds to considering the first returns at time $n$,
while $T_n$ considers all returns at time $n$. It is therefore
easy to check the following \emph{renewal equation} (see
e.g.~\cite[Proposition 1]{sarig:decay}):
  \begin{equation}
  \label{eq:renewal}
  T_n=\sum_{l=1}^\infty \sum_{k_1+\dots +k_l=n} R_{k_1}\dots
  R_{k_l}.
  \end{equation}
Let $\Lp$ be the space of functions of bounded variation on
$Y$. An element $f$ of $\Lp$ is a bounded function on $\R$,
supported in $Y$, and its norm is
  \begin{equation}
  \Var(f):=\sup_{N\in \N} \sup_{x_0<\dots<x_N} \sum_{i=0}^{N-1} |f(x_{i+1}) -f(x_i)|,
  \end{equation}
where the $x_i$'s are real numbers (not necessarily in $Y$). In
particular, $\norm{f}_{L^\infty}\leq \Var(f)/2$.

\begin{lem}
\label{lem:SatisfiesRenewal} The operators $R_n$ acting on $\Lp$
satisfy the assumptions of Theorem \ref{thm:renewal}. The spectral
projection $P$ corresponding to the eigenvalue $1$ of $R(1)$ is
given by
  \begin{equation}
  \label{eq:GivesP}
  Pf=\frac{ \left(\int_Y f\dLeb\right)}{\mu(Y)} h_Y
  \end{equation}
where $h_Y$ is the restriction to $Y$ of the density $h$ of the
invariant probability measure $\mu$. Additionally,
$PR'(1)P=P/\mu(Y)$.
\end{lem}
\begin{proof}
This lemma is proved in \cite{gouezel:decay} for the action of $R_n$
on the space $\TLp$ of Lipschitz functions on $Y$. We will adapt
this proof to the space $\Lp$.

The set $\{\phi=n\}$ is a subinterval $I_n$ of $Y$, and $T^n$ is a
diffeomorphism between $I_n$ and $Y$. Moreover, $|I_n|\sim
c/n^{\beta+1}$ for some constant $c>0$, and the distortion of $T^n$
on $I_n$ is uniformly bounded, independently of $n$, in the
following sense: there exists $C>0$ such that, for all $x,y\in I_n$,
  \begin{equation}
  \left| 1-\frac{(T^n)'(x)}{(T^n)'(y)}\right|\leq C|T^n x- T^n y|.
  \end{equation}
See e.g.~\cite[Section 6]{lsyoung:recurrence} for a proof of
these facts. Let $\psi_n: Y\to I_n$ be the inverse of $T^n$ on
$I_n$, so that
  \begin{equation}
  R_n f(x)= \psi_n'(x) f(\psi_n x).
  \end{equation}
Then
  \begin{equation}
  \Var(R_n f)\leq \norm{\psi_n'}_{L^\infty}\Var(f\circ \psi_n)+
  \norm{f}_{L^\infty} \Var(\psi_n') \leq C |I_n|\Var(f).
  \end{equation}
In particular,
  \begin{equation}
  \label{eq:BonnenormeBV}
  \norm{R_n}_{\Lp\to \Lp}\leq \frac{C}{n^{\beta+1}}.
  \end{equation}

As in Theorem \ref{thm:renewal}, we define for $z\in
\overline{\D}$ an operator $R(z)=\sum R_n z^n$. By
\eqref{eq:BonnenormeBV}, this operator is well defined on
$\Lp$. Moreover, by \cite[Paragraph 6.3]{gouezel:decay}, $R(z)$
also acts continuously on the space $\TLp$ of Lipschitz
continuous functions on $Y$, and satisfies the following
properties. First of all, $R(z)$ satisfies a Lasota-Yorke
inequality between $\TLp$ and $L^1$. Hence, by the theorem of
Ionescu-Tulcea and Marinescu, any eigenfunction of $R(z)$ (for
an eigenvalue of modulus $1$) which belongs to $L^1$ belongs in
fact to $\TLp$. Moreover, for $z\in \overline{\D}\moins \{1\}$,
$I-R(z)$ is invertible on $\TLp$, while $R(1)$ has a simple
eigenvalue at $1$, the corresponding eigenfunction being $h_Y$

Let us now prove that, for any $z\in \overline{\D}$, the
essential spectral radius of $R(z)$ acting on $\Lp$ is $<1$.
This could be proved by mimicking the arguments in
\cite{rychlik}, but it is easier to refer to \cite[Theorem
B.1]{ruelle:funct_zeta_milnor}. Indeed, this theorem shows that
the essential spectral radius of $R(z)$ is bounded by
$\norm{z^\phi/(T^\phi)'}_{L^\infty}<1$.

Let $z\in \overline\D \moins\{1\}$. If $I-R(z)$ were not
invertible on $\Lp$, then there would exist a function $f\in
\Lp$ such that $R(z)f=f$. The function $f$ would in particular
belong to $L^1$, hence, by the above argument, it would belong
to $\TLp$. This is a contradiction since $I-R(z)$ is invertible
on $\TLp$. In the same way, we check that $R(1)$ has a simple
eigenvalue at $1$, the eigenfunction still being the density of
the invariant measure. Moreover, the eigenprojection is given
by \eqref{eq:GivesP}.

We compute finally $PR'(1)P$. The formula for $Pf$ gives
  \begin{equation}
  PR'(1)Pf = \frac{\left(\int_Y R'(1)h_Y \dLeb\right)}{\mu(Y)}
  \frac{\left( \int_Y f\dLeb\right)}{\mu(Y)} h_Y
  =\gamma Pf,
  \end{equation}
for $\gamma=\left(\int_Y R'(1)h_Y \dLeb\right)/\mu(Y)$.
Moreover,
  \begin{equation}
  \int R_n h_Y \dLeb=\int \hat{T}^n( 1_{\{\phi=n\}}h_Y) \dLeb= \int
  1_{\{\phi=n\}} h_Y \dLeb= \mu\{\phi=n\}.
  \end{equation}
Summing these formulas over $n$ gives
  \begin{equation}
  \int R'(1) h_Y \dLeb=\sum n \mu\{\phi=n\}= \int_Y \phi \dd\mu=1
  \end{equation}
by Kac Formula. Hence, $\gamma=1/\mu(Y)$.
\end{proof}

\begin{cor}
\label{cor:GoodDecay} There exist $C>0$, and a sequence $c_n$
of complex numbers converging to $1$ when $n$ tends to
infinity, such that, for any functions $f,g$ supported in $Y$,
for any $n>0$,
  \begin{equation}
  \left| \int f \cdot g\circ T^n \dLeb - c_n \left(\int
  f\dLeb\right)\left(\int g \dd\mu\right)\right| \leq \frac{C
  \norm{f}_\Lp \norm{g}_{L^1(\Leb)}}{ n^\beta}.
  \end{equation}
\end{cor}
\begin{proof}
We have
  \begin{equation}
  \label{eq:Reeexprime}
  \int f \cdot g\circ T^n \dLeb= \int 1_Y \hat{T}^n(1_Y f) g
  \dLeb=\int T_nf \cdot g\dLeb.
  \end{equation}
Moreover, by \eqref{eq:renewal}, Lemma
\ref{lem:SatisfiesRenewal} and Theorem \ref{thm:renewal}, there
exist a sequence $c_n$ converging to $1$ and a constant $C$
such that
  \begin{align*}
  \norm{ T_n f - c_n \left(\int_Y f \dLeb\right) h_Y }_{\Lp}
  & = \norm{ T_n f - c_n \mu(Y) Pf }_{\Lp}
  \leq \norm{f}_\Lp \norm{T_n-c_n \mu(Y) P}
  \\&
  \leq \frac{C \norm{f}_\Lp}{n^\beta}.
  \end{align*}
Together with \eqref{eq:Reeexprime}, this concludes the proof.
\end{proof}

\begin{proof}[Proof of Theorem \ref{main_thm}]
Let first $A_n$ be a sequence of intervals contained in
$(1/2,1]$, with $\sum\Leb(A_n)=\infty$ (or, equivalently, $\sum
\mu(A_n)=\infty$). Let $B_n=T^{-n}A_n$. Let $j>i$. Applying
Corollary \ref{cor:GoodDecay} to $f=1_{A_i}h_Y$, $g=1_{A_j}$
and $n=j-i$, we get
  \begin{align*}
  \left| \mu(B_i \cap B_j) - c_{j-i} \mu(B_i)\mu(B_j)\right|
  \!\!\!\!\!\!\!\!\!\!\!\!\!\!\!\!\!\!\!\!\!\!\!\!\!\!\!\!\!\!\!
  \!\!\!\!\!\!\!\!\!\!\!\!\!\!\!\!\!\!\!\!\!\!\!\!&
  \\&
  = \left| \int 1_{A_i}h_Y\cdot 1_{A_j}\circ T^{j-i}\dLeb -c_{j-i}\left(
  \int 1_{A_i}h_Y \dLeb\right) \left( \int 1_{A_j}\dd\mu\right)
  \right|
  \\&
  \leq \frac{C \Var(1_{A_i}h_Y) \Leb(A_j)}{(j-i)^\beta}.
  \end{align*}
The function $h$ is Lipschitz continuous on $Y$, and bounded
from below. In particular, $\Leb(A_j)\leq C\mu(A_j)=C\mu(B_j)$.
We conclude
  \begin{equation}
  \left| \mu(B_i \cap B_j) - c_{j-i} \mu(B_i)\mu(B_j)\right|
  \leq \frac{C \mu(B_j)}{(j-i)^\beta}.
  \end{equation}

Let $\epsilon>0$. Let $K$ be such that, for $n\geq K$,
$|c_n|\leq 1+\epsilon$. Then
  \begin{multline*}
  \sum_{0\leq i < j<n} \mu(B_i \cap B_j)
  \leq \sum_{0\leq i<j<n} |c_{j-i}| \mu(B_i)\mu(B_j) +
  \sum_{j=1}^{n-1} \left( \sum_{i=0}^{j-1}
  \frac{C'}{(j-i)^\beta}\right) \mu(B_j)
  \\
  \leq \sum_{0\leq i<j<n} (1+\epsilon) \mu(B_i)\mu(B_j)
  + \sum_{j=1}^{n-1} \left( K \sup_{p\in \N}|c_p|+\sum_{p=1}^{\infty}
  \frac{C'}{p^\beta}\right) \mu(B_j).
  \end{multline*}
Therefore,
  \begin{equation}
  \frac{ \sum_{0\leq i < j<n} \mu(B_i \cap
  B_j)}{\left(\sum_{j=0}^{n-1} \mu(B_j)\right)^2} \leq
  \frac{1+\epsilon}{2} + \left( K \sup_{p\in \N}|c_p|+\sum_{p=1}^{\infty}
  \frac{C'}{p^\beta}\right) \frac{1}{\sum_{j=0}^{n-1} \mu(B_j)}.
  \end{equation}
Since $\sum_{j\in \N} \mu(B_j)=\infty$, this upper bound is at most
$1/2+\epsilon$ for large enough $n$. We have proved that
  \begin{equation}
  \limsup_{n\to\infty} \frac{ \sum_{0\leq i < j<n} \mu(B_i \cap
  B_j)}{\left(\sum_{j=0}^{n-1} \mu(B_j)\right)^2} \leq \frac{1}{2}.
  \end{equation}
By Theorem \ref{thm:BC_abstrait}, this concludes the proof in this
case.

Consider now $A_n$ an arbitrary sequence of intervals in $(0,1]$
with $\sum \Leb(A_n)=\infty$. Let $A'_n=T^{-1}(A_{n+1})\cap
(1/2,1]$. Since $\Leb(A'_n)=\Leb(A_n)/2$, this sequence of intervals
satisfies $\sum \Leb(A'_n)=\infty$, and $A'_n$ is a subinterval of
$(1/2,1]$. The first part of the proof shows that, for almost every
$x$, $T^n x$ belongs to $A'_n$ infinitely often. However, if
$T^n(x)\in A'_n$, then $T^{n+1}(x) \in A_{n+1}$. This concludes the
proof.

\end{proof}
\bibliography{biblio}
\bibliographystyle{alpha}

\end{document}